\documentclass[10pt]{article}
\usepackage{graphicx}

\newtheorem{proposition}{Proposition}

\title{On a weighted quasi-residual minimization strategy of the QMR method for solving complex symmetric shifted linear systems\thanks{This 
        work was partially supported by KAKENHI (Grant No. 18760063, 19560065).}}

\author{
T. Sogabe\thanks{Department of Computational Science and Engineering, Nagoya University, Furo-cho, Chikusa-ku,  Nagoya 464-8603, Japan ({\tt \{sogabe,zhang\}@na.cse.nagoya-u.ac.jp}).}
\and  T. Hoshi\thanks{(1) Department of Applied Mathematics and Physics, Tottori University, Japan (2) Core Research for Evolutional Science and Technology, Japan Science and Technology Agency (CREST-JST), 4-1-8 Honcho, Kawaguchi-shi, Saitama 332-0012, Japan ({\tt hoshi@damp.tottori-u.ac.jp}).} 
\and  S.-L. Zhang$^\dagger$ \hspace{-3mm}
\and T. Fujiwara\thanks{(1) Center for Research and Development of Higher Education, The University of Tokyo, Hongo, 7-3-1, Bunkyo-ku, Tokyo 113-8656, Japan  (2) Core Research for Evolutional Science and Technology, Japan Science and Technology Agency (CREST-JST), 4-1-8 Honcho, Kawaguchi-shi, Saitama 332-0012, Japan ({\tt fujiwara@coral.t.u-tokyo.ac.jp}).} }

\newcommand{\bb} {\mbox{\boldmath $b$}}

\newcommand{\eb} {\mbox{\boldmath $e$}}

\newcommand{\gb} {\mbox{\boldmath $g$}}
 
\newcommand{\pb} {\mbox{\boldmath $p$}}
\newcommand{\rb} {\mbox{\boldmath $r$}}

\newcommand{\vb} {\mbox{\boldmath $v$}}
\newcommand{\wb} {\mbox{\boldmath $w$}}
\newcommand{\xb} {\mbox{\boldmath $x$}}
\newcommand{\yb} {\mbox{\boldmath $y$}}
\newcommand{\zb} {\mbox{\boldmath $z$}}
\newcommand{\txt}  {\textrm}
\newcommand{\TT}   {\textrm{\scriptsize{T}}}
\newcommand{\HH}   {\textrm{\scriptsize{H}}}
\newcommand{\hs}   {\hspace{-7mm}}

\begin{document}

\maketitle

\begin{abstract}
We consider the solution of complex symmetric shifted linear systems. Such systems arise in large-scale electronic structure simulations and there is a strong need for the fast solution of the systems. With the aim of solving the systems efficiently, we consider a special case of the QMR method for non-Hermitian shifted linear systems and propose its weighted quasi-minimal residual approach. A numerical algorithm, referred to as  shifted QMR\_SYM($B$), is given by the choice of a particularly cost-effective weight.  Numerical examples are presented to show the performance of the shifted QMR\_SYM($B$) method.


\end{abstract}

\section{Introduction}
In this paper we consider the solution of complex symmetric shifted linear systems of the form
\begin{eqnarray}
  (A+\sigma_\ell I)\xb^{(\ell)}=\bb, \quad   \ell=1,2,\ldots,m,  \label{CS_SLS}
\end{eqnarray}
where  $A({\sigma_\ell}):=A+\sigma_\ell I$ are nonsingular $N$-by-$N$ complex symmetric sparse matrices, i.e., $A(\sigma_\ell)=A(\sigma_\ell)^\TT\neq \overline{A}(\sigma_\ell)^\TT$, with scalar shifts $\sigma_\ell \in \mathbf{C}$, $I$ is the $N$-by-$N$ identity matrix, and $\xb^{(\ell)}, \bb$ are complex vectors of length $N$. The above systems arise in large-scale electronic structure simulations \cite{Shifted_COCG}, and there is a strong need for the fast solution of the systems.

Since the given shifted linear systems (\ref{CS_SLS}) are a set of sparse linear systems, it is natural to use Krylov subspace methods, and moreover since the coefficient matrices are complex symmetric, one of the simplest ways to solve the shifted linear systems is applying (preconditioned) Krylov subspace methods for solving complex symmetric linear systems such as the COCG method \cite{COCG}, the COCR method \cite{COCR}, and the QMR\_SYM method \cite{QMR_SYM} to all of the shifted linear systems (\ref{CS_SLS}). On the other hand, denoting the $n$-dimensional Krylov subspace with respect to $A$ and $\bb$ as $K_n(A,\bb):=\txt{span}\{\bb,A\bb,\ldots,A^{n-1}\bb\}$, we observe that 
\begin{eqnarray}
  K_n(A,\bb)=K_n(A({\sigma_\ell}),\bb).       \label{KS}
\end{eqnarray}
This implies that once basis vectors are generated for one of the Krylov subspaces $K_n(A({\sigma_\ell}),\bb)$, these basis vectors can be used to solve all the shifted linear systems.  In other words, there is no need to generate all Krylov subspaces $K_n(A(\sigma_\ell),\bb)$, and thus computational costs involving the basis generation, e.g., matrix-vector multiplications, are saved. Here we give a concrete example: if we apply the conjugate orthogonal conjugate gradient (COCG) method to  all the linear systems (\ref{CS_SLS}), then bases for $K_n(A(\sigma_\ell),\bb)$  are generated for $\ell=1,2,\ldots,m$.  On the other hand, if we apply the COCG method to just one of the shifted linear systems (\ref{CS_SLS}) (referred to as the ``$seed\ system$''), then the Krylov basis vectors are generated from the seed system and these vectors are used to solve the rest of the shifted linear systems. 

Based on the observation (\ref{KS}), the shifted COCG method \cite{Shifted_COCG}  has been recently proposed for solving complex symmetric shifted linear systems. The feature of the shifted COCG method is that the method performs COCG on a seed system and makes it possible to complete COCG for all shifted linear systems without further matrix-vector multiplications. The feature is completely different from some of the other well-known shifted linear solvers such as the shifted BiCGStab($\ell$) method \cite{Shifted_BiCGSSTABL} or the restarted shifted GMRES method \cite{Shifted_GMRES} since these perform  BiCGStab($\ell$) (or GMRES) on a seed system but a different method on the shifted linear systems in order to keep the residuals colinear.  The feature of the shifted COCG method plays a very important role in the seed switching technique \cite{Seed_switch} that avoids a minor problem of the shifted COCG method:  one requires the choice of a seed system and  an unsuitable choice may lead to the undesirable result that some shifted linear systems remain unsolved.



There is another approach to solving the shifted linear systems (\ref{CS_SLS}). That is the use of Krylov subspace methods for non-Hermitian shifted linear systems such as the shifted BiCGStab($\ell$) method \cite{Shifted_BiCGSSTABL}, the shifted (TF)QMR method \cite{Freund}, the restarted shifted FOM method \cite{Shifted_FOM}, and the restarted shifted GMRES method \cite{Shifted_GMRES}, see also, e.g., \cite{Simoncini}. We readily see that the relation (\ref{KS}) holds not only for complex symmetric matrices but also for non-Hermitian matrices, and these methods are based on the use of this shift-invariant relation. Therefore, this can be a good approach. However, since these methods do not exploit the property of complex symmetric matrices, their computational costs can be more expensive than that of the shifted COCG method.    

In this paper we consider the shifted QMR\_SYM method that is a special case of the QMR method for non-Hermitian shifted linear systems \cite{Freund} and clarify the most time consuming part of it for a large number of shifted linear systems. Then, in order to reduce the cost, we propose a weighted quasi-minimal residual (WQMR) approach and propose a specific weight. We experimentally show the practical efficiency of the resulting algorithm, referred to as shifted QMR\_SYM($B$),  when the number of shifted linear systems is large enough.


The present paper is organized as follows: in the next section, shifted QMR\_SYM  is described in order to specify the most time consuming part for a large number of shifted linear systems. In Section 3,  we propose a WQMR approach with a specific weight for reducing the cost of the most time consuming part. The resulting algorithm, shifted QMR\_SYM($B$), and its properties are given. In Section 4, some results of numerical examples from electronic structure simulations are shown to see the performance of the shifted QMR\_SYM($B$) method. Finally, some concluding remarks are made in Section 5.

Throughout this paper, unless otherwise stated, all vectors and matrices are assumed to be complex. $\overline{M}$, $M^\TT$, $M^\HH=\overline{M}^\TT$ denote the complex conjugate, transpose, and Hermitian transpose matrix of the matrix $M$. $\|\vb\|_W$ denotes the $W$-norm written as $(\vb^\HH W \vb)^{1/2}$, where $W$ is Hermitian positive definite.
\section{The QMR\_SYM method for solving complex symmetric shifted linear systems}
In this section, the shifted QMR\_SYM method and its properties for solving complex symmetric shifted linear systems are introduced. 

The QMR method for shifted linear systems was first formulated in \cite{Freund} for the case of a general non-Hermitian matrix. Therefore, by simplifying the non-Hermitian Lanczos process \cite{Lanczos}, as is known from other papers such as \cite{QMR_SYM,COCG}, a shifted simplified QMR method, shifted QMR\_SYM, is readily obtained for the case of a complex symmetric matrix. Its algorithm is given next. \\[2mm]
\fbox{\textbf{Algorithm 2.1. Shifted QMR\_SYM}}\\[5mm]
\begin{eqnarray*} 
&&\hs                                                 \\[-15mm]
&&\hs \xb_{0}^{(\ell)}=\pb_{-1}^{(\ell)}=
      \pb_{0}^{(\ell)}=\mathbf{0}, \ 
      \vb_{1}=\bb/(\bb^\TT \bb)^{1/2},\ g_1^{(\ell)}=(\bb^\TT \bb)^{1/2},\\
&&\hs \txt{{\bf for}}\ n=1,2, \ldots\ \txt{{\bf do:}}\\
&&\hs \quad \txt{(The complex symmetric Lanczos process)}\nonumber\\
&&\hs \quad \alpha_n =\vb_n^\TT A\vb_n,                           \nonumber \\
&&\hs \quad \widetilde{\vb}_{n+1}=
        A\vb_n-\alpha_n\vb_n-\beta_{n-1}\vb_{n-1},              \nonumber \\
&&\hs \quad \beta_{n} =(\widetilde{\vb}_{n+1}^\TT
                    \widetilde{\vb}_{n+1})^{1/2},                 \nonumber \\
&&\hs \quad \vb_{n+1} = \widetilde{\vb}_{n+1}/\beta_{n},          \nonumber \\
&&\hs \quad t_{n-1,n}^{(\ell)}=\beta_{n-1},\ t_{n,n}^{(\ell)} =\alpha_{n}+\sigma_{\ell},\  t_{n+1,n}^{(\ell)}=\beta_{n},         \nonumber \\
&&\hs \quad \txt{(Solve least squares problems by Givens rotations)}                            \nonumber \\
&&\hs \quad  \txt{{\bf for}}\ \ell=1,2, \ldots, m\ \txt{{\bf do:}}\\
&&\hs \qquad \txt{{\bf if}}\ \|\rb_{n}^{(\ell)}\|_2/\|\bb\|_2> \epsilon,\
            \txt{{\bf then}}\\
&&\hs \qquad \quad \txt{{\bf for}}\ i=\max\{1,n-2\}, \ldots, n-1\ 
                   \txt{{\bf do:}}\\
&&\hs \qquad \qquad \pmatrix{
                t_{i,n}^{(\ell)} \cr
               t_{i+1,n}^{(\ell)}}
=
\pmatrix{
	       c_{i}^{(\ell)}     & s_{i}^{(\ell)}  \cr
	       -\overline{s}_{i}^{(\ell)} & c_{i}^{(\ell)} }
\pmatrix{
               t_{i,n}^{(\ell)} \cr
               t_{i+1,n}^{(\ell)} },      \nonumber\\
&&\hs \qquad \quad \txt{{\bf end}}             \nonumber\\
&&\hs \qquad \quad c_n^{(\ell)}=\frac{|t_{n,n}^{(\ell)}|}{\sqrt{|t_{n,n}^{(\ell)}|^2+|t_{n+1,n}^{(\ell)}|^2}}, \nonumber \\
&&\hs \qquad  \quad \overline{s}_n^{(\ell)}=\frac{t_{n+1,n}^{(\ell)}}{t_{n,n}^{(\ell)}}c_n^{(\ell)},   \nonumber \\	
&&\hs \qquad  \quad  t_{n,n}^{(\ell)}  =c_n^{(\ell)}t_{n,n}^{(\ell)}+s_n^{(\ell)}t_{n+1,n}^{(\ell)},  \nonumber \\
&&\hs \qquad  \quad t_{n+1,n}^{(\ell)}=0,                        \nonumber \\
&&\hs \qquad \quad \pmatrix{
                g_{n}^{(\ell)} \cr
                g_{n+1}^{(\ell)} }
 =
 \pmatrix{
	       c_{n}^{(\ell)}     & s_{n}^{(\ell)}  \cr
	       -\overline{s}_{n}^{(\ell)} & c_{n}^{(\ell)} }
\pmatrix{
               g_{n}^{(\ell)} \cr
               0 },   \nonumber\\	
&&\hs \qquad \quad  \txt{(Update approximate solutions $\xb_n^{(\ell)}$)}\nonumber\\
&&\hs \qquad \quad \pb_{n}^{(\ell)}= \vb_{n}-(t_{n-2,n}^{(\ell)}/t_{n-2,n-2}^{(\ell)})\pb_{n-2}^{(\ell)}- (t_{n-1,n}^{(\ell)}/t_{n-1,n-1}^{(\ell)})\pb_{n-1}^{(\ell)}, \nonumber\\
&&\hs \qquad \quad \xb_{n}^{(\ell)}=\xb_{n-1}^{(\ell)}+(g_n^{(\ell)}/t_{n,n}^{(\ell)})\pb_n^{(\ell)},  \nonumber \\
&&\hs \quad \quad \txt{{\bf end if}}       \nonumber	\\
&&\hs \quad \txt{{\bf end}}                \nonumber	\\&&\hs \qquad \txt{{\bf if}}\ \|\rb_{n}^{(\ell)}\|_2/\|\bb\|_2\le \epsilon \ 
             \txt{for all}\ \ell,\ \txt{{\bf then} exit}.\\
&&\hs \txt{{\bf end}}                      \nonumber		     
\end{eqnarray*}
Algorithm 2.1 can be regarded as a natural combination of the results given in \cite{QMR_SYM,Freund}.

In order to know that the numerical solution is accurate enough, one may need to compute the residual 2-norms. In that case, the following computation may be useful.
\begin{proposition}[See \cite{FrNa91}]
The 2-norms of the $n$th residuals  of the approximate solutions $\xb_n^{(\ell)}$ of the shifted QMR\_SYM method are given by 
\begin{eqnarray*}
 \|\rb_n^{(\ell)}\|_2=|g_{n+1}^{(\ell)}|\cdot\|\wb_{n+1}^{(\ell)}\|_2\quad \txt{for} \quad   \ell=1,2,\ldots,m, 
\end{eqnarray*}
where $\wb_{n+1}^{(\ell)}=-s_n^{(\ell)}\wb_n^{(\ell)}+c_n^{(\ell)}\vb_{n+1}$ and $\wb_1^{(\ell)}=\vb_1.$
\end{proposition}

Proposition 2.1 is a result known to hold for the QMR method \cite{FrNa91}. Therefore, it also holds for the above specialized variant.
The rest of this section describes some special properties of the shifted QMR\_SYM method.  

\begin{proposition}[See \cite{Fr90}]
Let $A\in \mathbf{R}^{N\times N}$ be real symmetric, $\sigma_\ell\in\mathbf{C}$ be complex shifts, and $\bb\in\mathbf{R}^N$. Then the shifted QMR\_SYM method (Algorithm 2.1) enjoys the following properties:
\begin{itemize}
 \item[I.] all matrix-vector multiplications can be done in real arithmetic;
 \item[I\hspace{-0.5mm}I.] an approximate solution at $n$th iteration step for each $\ell$ has minimal residual 2-norms, i.e., $\xb_n^{(\ell)}$'s are generated such that 
$\min\|\rb_n^{(\ell)}\|_2$ over $\xb_n^{(\ell)}\in K_n(A,\bb)$;
\item[I\hspace{-0.5mm}I\hspace{-0.5mm}I.] $ \|\rb_n^{(\ell)}\|_2=|g_{n+1}^{(\ell)}|\  \txt{for} \  \ell=1,2,\ldots,m, \  n\ge 0.$ 
\end{itemize}
\end{proposition}
The above properties are known results since the properties have been proved for each individual shift. See \cite{Fr90} for details. 
The properties of Proposition 2.2 may be very useful for large-scale electronic structure simulations \cite{Shifted_COCG} and a projection approach for eigenvalue problems \cite{Sakurai} since there are complex symmetric shifted linear systems to be solved efficiently under the assumption of Proposition 2.2.  

\section{An iterative method for solving  complex symmetric shifted linear systems}
In this section we consider complex symmetric shifted linear systems with a large number of shifts. For such systems, say $m\gg 1$, the most time consuming part of Algorithm 2.1 can be generating approximate solutions since the cost for the recurrences is $6Nm+3m$ per iteration step. In order to reduce the cost, in this section we propose a weighted quasi-minimal residual approach with a specific weight. In the next subsection we discuss the details of the approach. In Subsection 3.2 we give a specific weight to achieve the reduction of the cost.

\subsection{A weighted quasi-minimal residual approach}
Before the formulation of a Weighted Quasi-Minimal Residual (WQMR) approach, let us recall the complex symmetric Lanczos process (see, e.g., Algorithm 2.1 in \cite{QMR_SYM}). \\[2.5mm]
\fbox{\textbf{Algorithm 3.1. The complex symmetric Lanczos process }}\\[-5mm]
\begin{eqnarray}
&&\hs\txt{set}\ \beta_0=0,\ \vb_0=\mathbf{0},\ \rb_0\ne \mathbf{0}\in \mathbf{C}^N,\ \nonumber \\
&&\hs\txt{set}\ \vb_1=\rb_0/(\rb_0^\TT\rb_0)^{1/2},         \nonumber \\
&&\hs\txt{for $n=1,2,\ldots, m-1$ do:}		   	    \nonumber \\
&&\hs\quad \alpha_n =\vb_n^\TT A\vb_n,                      \nonumber \\
&&\hs\quad \widetilde{\vb}_{n+1}=
        A\vb_n-\alpha_n\vb_n-\beta_{n-1}\vb_{n-1},              \nonumber \\
&&\hs\quad \beta_{n} =(\widetilde{\vb}_{n+1}^\TT \widetilde{\vb}_{n+1})^{1/2},   \nonumber \\
&&\hs\quad \vb_{n+1} = \widetilde{\vb}_{n+1}/\beta_{n}.               \nonumber \\
&&\hs\txt{end}    				                 \nonumber
\end{eqnarray}
The matrix form of Algorithm 3.1 is known as follows: let $T_{n+1,n}$ and $T_n$ be the $(n+1)\times n$ and $n\times n$ tridiagonal matrices whose entries are recurrence coefficients of the complex symmetric Lanczos process, which are given by 
\begin{eqnarray*}
T_{n+1,n}:=
\pmatrix{
	\alpha_1& \beta_1 &        &                \cr
	\beta_1 & \alpha_2& \ddots &                \cr
 	        & \ddots  & \ddots & \beta_{n-1}    \cr
	        &         & \beta_{n-1} & \alpha_n  \cr
	        &         &        & \beta_{n}},    \quad 
T_{n}:=
\pmatrix{
	\alpha_1& \beta_1 &        &                \cr
	\beta_1 & \alpha_2& \ddots &                \cr
	        & \ddots  & \ddots & \beta_{n-1}    \cr
	        &         & \beta_{n-1} & \alpha_n },
\end{eqnarray*}
and let $V_n$ be the $N\times n$ matrix with the Lanczos vectors as columns, i.e., $V_n:=(\vb_1,\vb_2,\ldots,\vb_n)$. Then from Algorithm 3.1, it follows that 
\begin{eqnarray}
&&\hspace{4mm} AV_n = V_{n+1}T_{n+1,n}  
  = V_n T_n +\beta_{n}\vb_{n+1}\eb^{\TT}_n, \label{comp_sym_Lanczos_Mat}
\end{eqnarray}
where $\eb_n=(0,0,\ldots,1)^\TT\in \mathbf{R}^n$.

Now we are ready to describe the WQMR approach.
Let $\xb_n^{(\ell)}$ be approximate solutions at the $n$th iteration step for the systems (\ref{CS_SLS}), which are given by
\begin{equation}
\xb_n^{(\ell)} = V_n\yb_n^{(\ell)},\quad \ell=1,2,\ldots,m, \label{sol}
\end{equation}
where $\yb_n^{(\ell)} \in \mathbf{C}^n$. Then, from the definition of residual vectors $\rb_n^{(\ell)}:=\bb-(A+\sigma_\ell I)\xb_n^{(\ell)}$, the update formulas (\ref{sol}), and the matrix form of the complex symmetric Lanczos process (\ref{comp_sym_Lanczos_Mat}), we readily obtain
\begin{eqnarray}
\rb_n^{(\ell)} 
 = V_{n+1}\Bigl(g_1\eb_1- T_{n+1,n}^{(\ell)}\yb_n^{(\ell)}\Bigr),\ \ \txt{where}\ \  T_{n+1,n}^{(\ell)}:= T_{n+1,n}+ \sigma_\ell
            \pmatrix{
              I_n \cr
              \mathbf{0}^\TT  \cr 
            }.                            \label{QMR_residual}
\end{eqnarray}
Here, $\eb_1$ is the first unit vector written by $\eb_1$=$(1,0,\ldots,0)^\TT$ and $g_1=(\bb^\TT\bb)^{1/2}$. It is natural to determine $\yb_n^{(\ell)}$ such that all residual 2-norms $\|\rb_n^{(\ell)}\|_2$ are minimized. However, this choice for $\yb_n^{(\ell)}$ is impractical due to large computational costs. Hence, an alternative approach is given, i.e., the vectors $\yb_n^{(\ell)}$ are determined by solving the following weighted least squares problems:
\begin{eqnarray}
\yb_n^{(\ell)}  = \arg\min_{\zb_n^{(\ell)}\in \mathbf{C}^n}\ \Bigl\| g_1\eb_1- T_{n+1,n}^{(\ell)}\zb_n^{(\ell)}\Bigr\|_{W_{n+1}^\HH W_{n+1}^{}},  
                                           \label{WLSP}
\end{eqnarray}
where $W_{n+1}$ is an $(n+1)$-by-$(n+1)$ nonsingular matrix. Thus $W_{n+1}^\HH W_{n+1}$ can be used as a weight since it is Hermitian positive definite.  
One of the simplest choices for $W_{n+1}$ is the identity matrix. In this case, from $W_{n+1}=I_{n+1}$ we have  
\begin{eqnarray}
\yb_n^{(\ell)}  = \arg\min_{\zb_n^{(\ell)}\in \mathbf{C}^n}\ \Bigl\| g_1 \eb_1- T_{n+1,n}^{(\ell)}\zb_n^{(\ell)}\Bigr\|_2.  \label{WLSP_2norm}
\end{eqnarray}
The vector that is minimized is called quasi-residual.
 Algorithm 2.1 is obtained by solving (\ref{WLSP_2norm}) using Givens rotations, see, e.g., \cite[p.215]{Mat_Comp}.

A slightly generalized choice proposed in \cite{FrNa91} is $W_{n+1}=\Omega_{n+1}:=\txt{diag}(\omega_1,\omega_2,\ldots,\\ \omega_{n+1})$ with $\omega_i>0$ for all $i$. Then, it follows that 
\begin{eqnarray*}
\yb_n^{(\ell)}  = \arg \min_{\zb_n^{(\ell)}\in \mathbf{C}^n}\ \Bigl\| \omega_1g_1\eb_1- \Omega_{n+1} T_{n+1,n}^{(\ell)}\zb_n^{(\ell)}\Bigr\|_2.
\end{eqnarray*}
Of various possible choices for $\omega_i$, a natural one is $\omega_i=\|\vb_i\|_2$ since $\Omega_{n+1}$ contains then the diagonal entries of the upper triangular matrix $R_{n+1}$ that is obtained by the $QR$ factorization of $V_{n+1}$. If we choose $W_{n+1}=R_{n+1}$, where $V_{n+1}=Q_{n+1}R_{n+1}$, then from (\ref{QMR_residual}) and (\ref{WLSP}) we have 
\begin{eqnarray*}
\min_{\zb_n^{(\ell)}\in \mathbf{C}^n}\ \Bigl\|g_1\eb_1- T_{n+1,n}^{(\ell)}\zb_n^{(\ell)} \Bigr\|_{R_{n+1}^\HH R_{n+1}^{}}
&=&\min_{\zb_n^{(\ell)}\in \mathbf{C}^n}\ \Bigl\|g_1R_{n+1}\eb_1- R_{n+1}T_{n+1,n}^{(\ell)}\zb_n^{(\ell)} \Bigr\|_2 \\
&=&\min_{\zb_n^{(\ell)}\in \mathbf{C}^n}\ \Bigl\|Q_{n+1}R_{n+1}\bigl(g_1\eb_1- T_{n+1,n}^{(\ell)}\zb_n^{(\ell)}\bigr)\Bigr\|_2 \\
&=&\min_{\zb_n^{(\ell)}\in \mathbf{C}^n}\ \Bigl\|V_{n+1}\bigl(g_1\eb_1- T_{n+1,n}^{(\ell)}\zb_n^{(\ell)}\bigr)\Bigr\|_2 \\
&=&\min_{\zb_n^{(\ell)}\in \mathbf{C}^n}\ \bigl\|\rb_n^{(\ell)}\bigr\|_2. \\
\end{eqnarray*}
By solving the above weighted least squares problems, all residual 2-norms are minimized. Hence  $W_{n+1}=\Omega_{n+1}$ is a rational choice. For each individual shift, the resulting algorithm is the same as that in \cite[Algorithm 3.2]{QMR_SYM}.
\subsection{A choice of the weight suitable for a large number of shifts}
In the previous subsection, we have described the WQMR approach and mentioned that the choice of the weight $W_{n+1}^\HH W_{n+1}$  with  $W_{n+1}=I_{n+1}$ leads to the shifted QMR\_SYM method (Algorithm 2.1). Under the assumption of proposition 2.2, the shifted QMR\_SYM method is ideal in the sense of Faber-Manteuffel's theorem \cite{Faber} since it enjoys minimal residual property and requires not long-term but short-term recurrences for updating approximate solutions, and thus one may think that there is no need to choose other possible weights. However, we will show in this subsection that even under the assumption of proposition 2.2 there is a practical weight for the WQMR approach. The motivation for the choice of the weight mainly comes from the freedom of the number $m$ of complex symmetric shifted linear systems.

Now we consider the computational costs of Algorithm 2.1 for a large number of complex symmetric shifted linear systems, i.e., $m\gg1$. For the case $m\gg1$ we readily see from Algorithm 2.1 that computing the recurrences for updating approximate solutions is the most time-consuming part due to a cost of  $6Nm+3m$ per iteration step. Hence we will now consider a weight to reduce the computational cost for the recurrences of $\xb_n^{(\ell)}$. To achieve this we propose the following choice:
\begin{eqnarray}
  W_{n+1}=L_{n+1}^{(\ell)}\ \txt{such that}\ L_{n+1}^{(\ell)}T_{n+1,n}^{(\ell)}=
\pmatrix{B_n^{(\ell)} \cr \mathbf{0}^T},   \label{bidiagonal_weight}
\end{eqnarray} 
where $B_n^{(\ell)}$ is an $n$-by-$n$ upper bidiagonal matrix of the form
\begin{eqnarray*}
B_{n}^{(\ell)}:=
\pmatrix{
t_{1,1}^{(\ell)}& t_{1,2}^{(\ell)} &        &       \cr
                & t_{2,2}^{(\ell)} & \ddots &       \cr
                &          &         \ddots &   t_{n-1,n}^{(\ell)} \cr
                &          &                &   t_{n,n}^{(\ell)}},
\end{eqnarray*}
and $L^{(\ell)}_{n+1}$ is lower triangular and will be specified below.

Next we derive recurrence formulas for updating the approximate solutions $\xb_n^{(\ell)}$. From (\ref{WLSP}) with the choice $W_{n+1}=L_{n+1}^{(\ell)}$ of (\ref{bidiagonal_weight}) it follows that
\begin{eqnarray}
\yb_n^{(\ell)}  &=& \arg\min_{\zb_n^{(\ell)}\in \mathbf{C}^n}\ \Bigl\| g_1 \eb_1- T_{n+1,n}^{(\ell)}\zb_n^{(\ell)}\Bigr\|_{(L_{n+1}^{(\ell)})^\HH L_{n+1}^{(\ell)}} \label{WLSP_B}\\
                &=& \arg\min_{\zb_n^{(\ell)}\in \mathbf{C}^n}\ \Bigl\| g_1  L_{n+1}^{(\ell)}\eb_1-  L_{n+1}^{(\ell)}T_{n+1,n}^{(\ell)}\zb_n^{(\ell)}\Bigr\|_{2}      \nonumber\\
                &=& \arg\min_{\zb_n^{(\ell)}\in \mathbf{C}^n}\ \biggl\| 
\pmatrix{\widetilde{\gb}_n^{(\ell)} \cr \widetilde{g}_{n+1}^{(\ell)}}
- \pmatrix{B_n^{(\ell)} \cr \mathbf{0}^\TT}
  \zb_n^{(\ell)}\biggr\|_2, \ \pmatrix{\widetilde{\gb}_n^{(\ell)} \cr \widetilde{g}_{n+1}^{(\ell)}}:=g_1 L_{n+1}^{(\ell)}\eb_1. \nonumber
\end{eqnarray}
From the above least squares problems we readily see that $\yb_n^{(\ell)}=\ \bigr(B_n^{(\ell)}\bigl)^{-1}\widetilde{\gb}_n^{(\ell)}$. Hence it follows from (\ref{sol}) and using $(\widetilde{\pb}_1\ \widetilde{\pb}_2\ \cdots\ \widetilde{\pb}_n):=V_n(B_n^{(\ell)})^{-1}$ that we have the following coupled two-term recurrence relations:
\begin{eqnarray*}
 \widetilde{\pb}_{n}^{(\ell)}&=& (\vb_{n}-t_{n-1,n}^{(\ell)}\widetilde{\pb}_{n-1}^{(\ell)})/t_{n,n}^{(\ell)},   \label{x_rec1_B}\\
 \xb_{n}^{(\ell)}&=&\xb_{n-1}^{(\ell)}+\widetilde{g}_n^{(\ell)}\widetilde{\pb}_n^{(\ell)}. \label{x_rec2_B}
\end{eqnarray*}
The cost per iteration is now $5Nm$.
Substituting  $\pb_i^{(\ell)}=t_{i,i}^{(\ell)}\widetilde{\pb}_i^{(\ell)}$ into the above recurrences, we have the even more efficient recurrence formulas.
\begin{eqnarray*}
 \pb_{n}^{(\ell)}&=& \vb_{n}-(t_{n-1,n}^{(\ell)}/t_{n-1,n-1}^{(\ell)})\pb_{n-1}^{(\ell)}, \\
 \xb_{n}^{(\ell)}&=&\xb_{n-1}^{(\ell)}+(\widetilde{g}_n^{(\ell)}/t_{n,n}^{(\ell)})\pb_n^{(\ell)}.
\end{eqnarray*}
By this reformulation, the cost becomes $4Nm+2m$.

Next, let us consider a choice for $L_{n+1}^{(\ell)}$  satisfying (\ref{bidiagonal_weight}). Let $F_{n}^{(\ell)}$ be an $n$-by-$n$ matrix for the form
\begin{eqnarray}
  F_n^{(\ell)}(i):= 
  \pmatrix{
    I_{i-1}  &           &     & \cr
             &      1    &     & \cr
             &  f_{i}^{(\ell)}&   1 & \cr
             &           &     & I_{n-i-1}        
  },    \label{F}
\end{eqnarray}
and let $T$ be tridiagonal.
Then, by determining $f_i^{(\ell)}$ such that the $(i+1,i)$ entry of the matrix $F_n^{(\ell)}(i)T$ is zero, we can fulfill (\ref{bidiagonal_weight}) in the following way:
\begin{eqnarray}
  F_{n+1}^{(\ell)}(n)F_{n+1}^{(\ell)}(n-1)\cdots F_{n+1}^{(\ell)}(1)T_{n+1,n}^{(\ell)}=
\pmatrix{
         B_n^{(\ell)} \cr
         \mathbf{0}^\TT  \cr 
         }  \label{FN}.      
\end{eqnarray}
From the above we see that $F_{n+1}^{(\ell)}(n)F_{n+1}^{(\ell)}(n-1)\cdots F_{n+1}^{(\ell)}(1)=L_{n+1}^{(\ell)}$, and thus  $L_{n+1}^{(\ell)}$ and $L_{n}^{(\ell)}$ are related by
\begin{eqnarray}
  L_{n+1}^{(\ell)} = F_{n+1}^{(\ell)}(n)\pmatrix{
         L_n^{(\ell)} & \mathbf{0}\cr  
         \mathbf{0}^\TT  & 1\cr 
         } \    \txt{for}\  n=2,3,\ldots,  \label{Relation_L}
\end{eqnarray}
where $L_2^{(\ell)}=F_2^{(\ell)}(1).$ From the above we readily see that the matrices $L_{n+1}^{(\ell)}$ are lower triangular with all 1's on the diagonals,  and this property will be used for the proof of proposition 3.1 given later. 
The shifted QMR\_SYM method with the weight $(L_{n+1}^{(\ell)})^\HH L_{n+1}^{(\ell)}$ is referred to as shifted QMR\_SYM($B$). Its algorithm is given next.\\ \\ 
\fbox{\textbf{Algorithm 3.2. Shifted QMR\_SYM$(B)$}}\\[7mm]
\begin{eqnarray*} 
&&\hs                                                 \\[-15mm]
&&\hs \xb_{0}^{(\ell)}=\pb_{-1}^{(\ell)}=
      \pb_{0}^{(\ell)}=\mathbf{0}, \ 
      \vb_{1}=\bb/(\bb^\TT\bb)^{1/2},\ \widetilde{g}_1^{(\ell)}=(\bb^\TT\bb)^{1/2},\\
&&\hs \txt{{\bf for}}\ n=1,2, \ldots\ \txt{{\bf do:}}\\
&&\hs \quad \txt{(The complex symmetric Lanczos process)}\nonumber\\
&&\hs \quad \alpha_n =\vb_n^\TT A\vb_n,                           \nonumber \\
&&\hs \quad \widetilde{\vb}_{n+1}=
        A\vb_n-\alpha_n\vb_n-\beta_{n-1}\vb_{n-1},              \nonumber \\
&&\hs \quad \beta_{n} =(\widetilde{\vb}_{n+1}^\TT
                    \widetilde{\vb}_{n+1})^{1/2},                  \nonumber \\
&&\hs \quad \vb_{n+1} = \widetilde{\vb}_{n+1}/\beta_{n},          \nonumber \\
&&\hs \quad t_{n-1,n}^{(\ell)}=\beta_{n-1},\ t_{n,n}^{(\ell)} =\alpha_{n}+\sigma_{\ell},\  t_{n+1,n}^{(\ell)}=\beta_{n},         \nonumber \\
&&\hs \quad \txt{(Solve weighted least squares problems)}                            \nonumber \\
&&\hs \quad  \txt{{\bf for}}\ \ell=1,2, \ldots, m\ \txt{{\bf do:}}\\
&&\hs \qquad \txt{{\bf if}}\ \|\rb_{n}^{(\ell)}\|_2/\|\bb\|_2> \epsilon,\
            \txt{{\bf then}}\\
&&\hs \qquad \quad \txt{{\bf for}}\ i=\max\{1,n-1\}, \ldots, n-1\ 
                   \txt{{\bf do:}}\\
&&\hs \qquad \qquad t_{i+1,n}^{(\ell)} =  f_{i}^{(\ell)} t_{i,n}^{(\ell)} +   t_{i+1,n}^{(\ell)},   \nonumber\\
&&\hs \qquad \quad \txt{{\bf end}}             \nonumber\\
&&\hs \qquad  \quad f_n^{(\ell)}=-\frac{t_{n+1,n}^{(\ell)}}{t_{n,n}^{(\ell)}},   \nonumber \\	
&&\hs \qquad  \quad t_{n+1,n}^{(\ell)}=0,                        \nonumber \\
&&\hs \qquad \quad  \widetilde{g}_{n+1}^{(\ell)} =f_{n}^{(\ell)}\widetilde{g}_{n}^{(\ell)},  \nonumber\\
&&\hs \qquad \quad  \txt{(Update approximate solutions $\xb_n^{(\ell)}$)}\nonumber\\
&&\hs \qquad \quad \pb_{n}^{(\ell)}= \vb_{n}-(t_{n-1,n}^{(\ell)}/t_{n-1,n-1}^{(\ell)})\pb_{n-1}^{(\ell)}, \nonumber\\
&&\hs \qquad \quad \xb_{n}^{(\ell)}=\xb_{n-1}^{(\ell)}+(\widetilde{g}_n^{(\ell)}/t_{n,n}^{(\ell)})\pb_n^{(\ell)},  \nonumber \\
&&\hs \quad \quad \txt{{\bf end if}}       \nonumber	\\
&&\hs \quad \txt{{\bf end}}                \nonumber	\\&&\hs \qquad \txt{{\bf if}}\ \|\rb_{n}^{(\ell)}\|_2 / \|\bb\|_2\le \epsilon \ 
             \txt{for all}\ \ell,\ \txt{{\bf then} exit}.\\
&&\hs \txt{{\bf end}}                      \nonumber		     
\end{eqnarray*}
Similar to Proposition 2.1, there is an efficient way to evaluate residual 2-norms as follows:
\begin{proposition}
The $n$th residual 2-norms of the approximate solutions $\xb_n^{(\ell)}$ for the shifted QMR\_SYM($B$) method are given by 
\begin{eqnarray*}
 \|\rb_n^{(\ell)}\|_2=|\widetilde{g}_{n+1}^{(\ell)}|\cdot\|\vb_{n+1}\|_2\quad \txt{for} \quad   \ell=1,2,\ldots,m.
\end{eqnarray*}
\end{proposition}

{\em Proof}. The proof is similar to that of Proposition 2.1. It follows from (\ref{QMR_residual}), (\ref{bidiagonal_weight}), (\ref{WLSP_B}), and recalling  $\yb_n^{(\ell)}=\ \bigr(B_n^{(\ell)}\bigl)^{-1}\widetilde{\gb}_n^{(\ell)}$ that we have
\begin{eqnarray*}
\rb_{n}^{(\ell)}=\widetilde{g}_{n+1}^{(\ell)}V_{n+1}\bigr(L_{n+1}^{(\ell)}\bigl)^{-1}\eb_{n+1}. 
\end{eqnarray*}
From (\ref{F}) and (\ref{FN}) $L_{n+1}^{(\ell)}$ is a lower triangular matrix with all 1's on the diagonals. Thus, $\bigr(L_{n+1}^{(\ell)}\bigl)^{-1}$ is also a lower triangular matrix with all 1's on the diagonals. It follows that  $\bigr(L_{n+1}^{(\ell)}\bigl)^{-1}\eb_{n+1}=\eb_{n+1}$, and thus we have $\rb_{n}^{(\ell)}=\widetilde{g}_{n+1}^{(\ell)}V_{n+1}\eb_{n+1}=\widetilde{g}_{n+1}^{(\ell)}\vb_{n+1}$, which concludes the proof. 
\qquad 

When we solve  a large number of shifted systems, the computational cost of the residual 2-norms for the shifted QMR\_SYM($B$) method is much cheaper than the that for the shifted QMR\_SYM method since the former cost is of order $N$ and the latter is of order $Nm$ per iteration step.

Observing algorithms of the shifted QMR\_SYM and the shifted QMR\_SYM($B$) methods, we see that the work for the weighted least squares problems and for updating approximate solutions in the shifted QMR\_SYM($B$) method is lower than that in the shifted QMR\_SYM method. For the case $m\gg 1$, the most time-consuming part is updating the approximate solutions. In this case, the shifted QMR\_SYM($B$) method may be more efficient than the shifted QMR\_SYM method since the shifted QMR\_SYM($B$) method requires $4Nm+2m$ operations per iteration step for the update while the shifted QMR\_SYM method needs $6Nm+3m$ operations. This difference will be clearer when we use the results of Propositions 2.1 and 3.1 as stopping criterion. On the other hand, in terms of number of iterations, under a certain assumption, the convergence of the shifted QMR\_SYM($B$) method  is worse than that of the shifted QMR\_SYM method but not worse than that of  the shifted COCG method, as described by the following proposition:
\begin{proposition}
Under the assumption of Proposition 2.2, the shifted QMR\_SY\\M($B$) method (Algorithm 3.2) enjoys the following properties:
\begin{itemize}
 \item[I.] all matrix-vector multiplications can be done in real arithmetic;
 \item[I\hspace{-0.5mm}I.] if breakdown does not occur and each matrix $T_n+\sigma_\ell I_n$is nonsingular, then $\bigr\|\rb_n^{(\ell), \txt{\scriptsize{SQ($B$)}}}\bigl\|_2=\bigr\|\rb_n^{(\ell), \txt{\scriptsize{SCOCG}}}\bigl\|_2\ge \bigr\|\rb_n^{(\ell), \txt{\scriptsize{SQ}}}\bigl\|_2 \ \txt{for}\   \ell=1,2,\ldots m$, where the superscripts SQ($B$), SCOCG, and SQ are short for shifted QMR\_SYM($B$), shifted COCG, and shifted QMR\_SYM, respectively;
\item[I\hspace{-0.5mm}I\hspace{-0.5mm}I.] $ \|\rb_n^{(\ell), \txt{\scriptsize{SQ($B$)}}}\|_2=|\widetilde{g}_{n+1}^{(\ell)}|\  \txt{for} \  \ell=1,2,\ldots,m, \  n\ge 0.$ 
\end{itemize}
\end{proposition}
{\em Proof}.
The proof of ($I$) is the same as that of Proposition 2.2, and is based on the fact that under the assumption the complex symmetric Lanczos process generates real basis vectors.

Next, we give a proof of ($I\hspace{-0.5mm}I$). The $n$th residuals of the shifted COCG method \cite{Shifted_COCG} belong to $\bb-(A+\sigma_\ell I)K_n(A+\sigma_\ell I,\bb)$. Hence, each $\rb_n^{(\ell), \txt{\scriptsize{SCOCG}}}$ can be written as $\rb_n^{(\ell), \txt{\scriptsize{SCOCG}}}=\bb-(A+\sigma_\ell I)V_n\yb_n^{(\ell), \txt{\scriptsize{SCOCG}}}$, where $V_n$ is the same matrix as in (\ref{comp_sym_Lanczos_Mat}) Since each $\rb_n^{(\ell), \txt{\scriptsize{SCOCG}}}$ is orthogonal to each subspace $K_n(\overline{A}+\overline{\sigma}_\ell I,\overline{\bb})$\{=$K_n(\overline{A},\overline{\bb})\}$, i.e., $\rb_n^{(\ell), \txt{\scriptsize{SCOCG}}} \perp K_n(\overline{A},\overline{\bb})$, we have $V_n^\TT\bb-V_n^\TT(A+\sigma_\ell I)V_n\yb^{(\ell), \txt{\scriptsize{SCOCG}}}=0$, and thus it follows from (\ref{comp_sym_Lanczos_Mat}) that we have the relation $\yb_n^{(\ell), \txt{\scriptsize{SCOCG}}} = \{V_n^\TT(A+\sigma_\ell I)V_n\}^{-1}V_n^\TT\bb = g_1(T_n^{(\ell)})^{-1}\eb_1$, where $T_n^{(\ell)}:=T_n+\sigma_\ell I_n$. Since the shifted QMR\_SYM($B$) method has the form (\ref{QMR_residual}), it is sufficient to show that  $\yb_n^{(\ell), \txt{\scriptsize{SCOCG}}}=\yb_n^{(\ell), \txt{\scriptsize{SQ($B$)}}}$. From (\ref{WLSP_B}) and (\ref{Relation_L}) it follows that $\yb_n^{(\ell), \txt{\scriptsize{SQ($B$)}}}=(B_n^{(\ell)})^{-1}\widetilde{\gb}_n^{(\ell)}=g_1(B_n^{(\ell)})^{-1}[I_n |\mathbf{0}]L_{n+1}^{(\ell)}[\eb_1^\TT | 0]^\TT=g_1(B_n^{(\ell)})^{-1}L_{n}^{(\ell)}\eb_1=g_1\{(L_{n}^{(\ell)})^{-1}B_n^{(\ell)}\}^{-1} \eb_1.$ Since from (\ref{FN}) and (\ref{Relation_L}) we can readily confirm the relation $L_{n}^{(\ell)}T_n^{(\ell)}=B_n^{(\ell)}$, we have  $\yb_n^{(\ell), \txt{\scriptsize{SQ($B$)}}}=g_1(T_n^{(\ell)})^{-1} \eb_1,$ which is the same as $\yb_n^{(\ell), \txt{\scriptsize{SCOCG}}}.$ The inequality in ($I\hspace{-0.5mm}I$) follows from Proposition 2.2 since under the given assumption the shifted QMR\_SYM method enjoys the minimal residual property.

Finally, we give a proof of ($I\hspace{-0.7mm}I\hspace{-0.7mm}I$). If follows from the proof of ($I$) that  $\|\vb_i\|_2=1$ for all $i$. Thus from Proposition 3.1 we have $ \|\rb_n^{(\ell), \txt{\scriptsize{SQ$(B)$}}}\|_2=|\widetilde{g}_{n+1}^{(\ell)}|\cdot\|\vb_{n+1}\|_2=|\widetilde{g}_{n+1}^{(\ell)}|\ \txt{for} \   \ell=1,2,\ldots,m,\ n\ge 0.$
\qquad 

In the property ($I\hspace{-0.7mm}I$) of Proposition 3.2, breakdown may occur due to the choice (\ref{F}) of the weighted least squares problems.

From proposition 3.2 we see that in terms of the number of iteration steps the shifted QMR\_SYM($B$) method never converges faster than the shifted QMR\_SYM method but it converges at the same iteration step as the shifted COCG method does. Since the efficiency of the shifted COCG method has been shown already and for the case $m\gg 1$ the computational cost of the shifted QMR\_SYM$(B)$ method is much less than that of the shifted QMR\_SYM method, the shifted QMR\_SYM$(B)$ method can also be useful. This will be supported by some numerical examples in the next section.  
\section{Numerical examples}
In this section, we report on some numerical examples with the shifted COCG method, the shifted QMR\_SYM method (Algorithm 2.1), and the shifted QMR\_SYM($B$) method (Algorithm 3.2).  We evaluate these methods in terms of computational time.  
All tests were performed on a workstation with a 2.6GHz AMD Opteron(tm) processor 252 using double precision arithmetic. 
Codes were written in Fortran 77 and compiled with g77 -$O$3. 
In all cases the stopping 
criterion was set as $\epsilon=10^{-12}$.   \\ \\
3.1. Example 1\\ \\
The first problem comes from the electronic structure calculation of a bulk Si(001) with 512 atoms in \cite{Shifted_COCG}, which is written as follows: 
\begin{eqnarray*}
(\sigma_\ell I - H)\xb^{(\ell)} = \eb_1, \ \  \ell=1,2, \ldots,m,
\end{eqnarray*}
where $\sigma_\ell=0.400 + (\ell-1+i)/1000$, $H\in R^{2048\times 2048}$ is a symmetric matrix with 139264 entries,  $\eb_1=(1,0,\ldots,0)^\TT$, and $m=1001$. Since the shifted COCG method requires the choice of a seed system, we have chosen the optimal seed ($\ell=714$) in this problem; otherwise some linear systems will remain unsolved.


Figure 4.1 shows the number of iterations of each method to solve the $\ell$th shifted linear systems. For example, in Fig.~4.1, the number of iterations for the shifted COCG method at $\ell=600$ is 150, which means the shifted COCG method required 150 iterations to obtain the (approximate) solutions of the 600th shifted linear system, i.e., $(\sigma_{600} I - H)\xb^{(600)} = \eb_1$. 

\begin{figure}[htbp]
\begin{center}
\includegraphics[width=0.8\linewidth, height=0.5\linewidth]{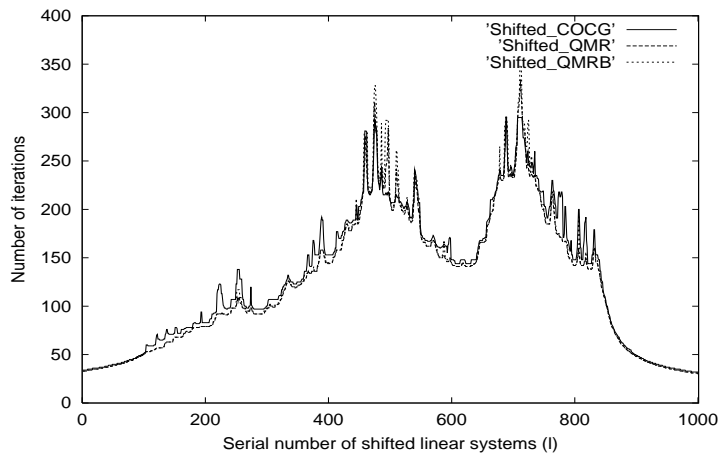}
\end{center}
\caption{Number of iterations for the shifted COCG method, the shifted QMR\_SYM method, and the shifted QMR\_SYM($B$) method versus the serial number of the shifted linear systems. } 
\end{figure}

\begin{figure}[htbp] 
\begin{center}
\includegraphics[width=0.8\linewidth, height=0.5\linewidth]{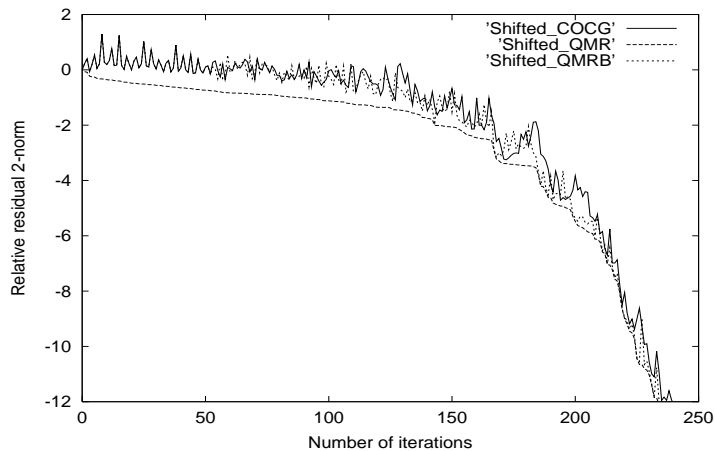}
\end{center}
\caption{{\rm Log}$_{10}$ of the relative residual 2-norms versus the number of iterations of the shifted COCG method, the shifted QMR\_SYM method, and the shifted QMR\_SYM($B$) method for the shifted linear system with $\ell=701, i.e., \sigma_{701}=1.100+0.001i$.} 
\end{figure} 

From Fig.~4.1 we make three observations: first, the three methods required almost the same number of iterations at each $\ell$; second, in terms of number of iterations, the shifted QMR\_SYM method often converged slightly faster than the other two methods. This phenomenon is closely related to Proposition 2.2, which will become clearer later; third, for each method the required number of iterations depends highly on the shift parameters $\sigma_\ell$. This result may come from varying eigenvalues of the coefficient matrices $\sigma_\ell I -H$  since if we choose $\sigma_\ell$ close to an eigenvalue of $H$, then $\sigma_\ell I -H$ is close to singular. Conversely, from the shape of Fig.~4.1 one may obtain a partial distribution of eigenvalues of $H$.

One of the residual 2-norm histories for the three methods is given in Fig.~4.2. From Fig.~4.2 we see that the relative residual 2-norm of the shifted QMR\_SYM method decreases monotonically and at every iteration step the norm is less than those of the other two methods. Hence we can say that the property $(I\hspace{-0.5mm}I)$ of Proposition 2.2 is experimentally supported by this history. We make another observation in the histories of the shifted COCG method and the shifted QMR\_SYM($B$) method: During the first about fifty iterations, the two methods show the same histories. After that their histories varies gradually. Hence we see that the property $(I\hspace{-0.5mm}I)$ of Proposition 3.2 is also experimentally supported by these histories.

\begin{figure}[htbp] 
\begin{center}
\includegraphics[width=0.8\linewidth, height=0.5\linewidth]{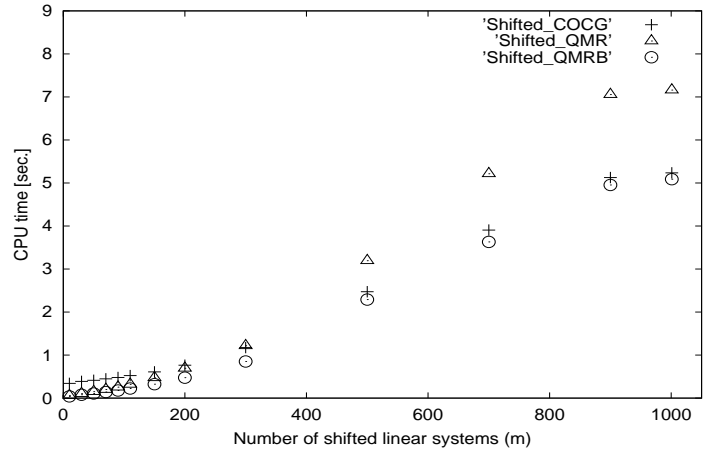}
\caption{Required CPU time given in seconds versus the number of shifted linear systems for each iterative method.} 
\end{center}
\end{figure} 

\begin{figure}[htbp] 
\begin{center}
\includegraphics[width=0.8\linewidth, height=0.5\linewidth]{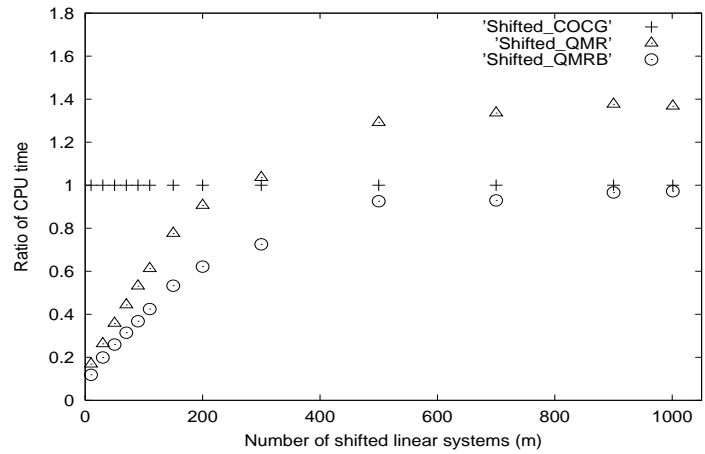}
\end{center}
\caption{The ratio of each computational time to the one of the shifted COCG method versus the number of shifted linear systems.} 
\end{figure} 

Each computational time of the three methods is given in Fig.~4.3, where the horizontal axis denotes the number $m$ of shifted linear systems that are solved. For example, in Fig.~4.3, the computational time of the shifted COCG method at $m=200$ is about 0.76 [sec.], which means that it required about 0.76 [sec.] to solve the shifted linear systems: $\{(0.400+0.001i) I - H\}\xb^{(1)} = \eb_1, \{(0.401+0.001i) I - H\}\xb^{(2)} = \eb_1, \ldots, \{(0.599+0.001i) I - H\}\xb^{(200)} = \eb_1$.  From Fig.~4.3 we see that as the number $m$ grows larger, the shifted QMR\_SYM method required more CPU time than the other two methods. On the other hand, the shifted QMR\_SYM($B$) method required almost the same CPU time as the shifted COCG method. This phenomena can be attributed to the computational costs of updating approximate solutions for each method since there are three facts: first, we know from Fig.~4.1 that the three methods required almost the same iterations; second, the shifted QMR\_SYM($B$) method requires almost the same computational costs as the shifted COCG method, while the shifted QMR\_SYM method tends to require more computational costs per iteration than the other two methods; third, for large $m$, updating approximate solutions is one of the most time-consuming parts. The second and the third facts have been discussed already in the previous two sections.

 In Fig.~4.3 we can see little about the properties of the three methods for small $\ell$.  We therefore show the ratio of each computational time to the computational time of the shifted COCG method in Fig.~4.4.
We see from Fig.~4.4 that in terms of the ratio of CPU times the shifted QMR\_SYM method and the shifted QMR\_SYM($B$) method converged much faster than the shifted COCG method when the number of shifted linear systems is small, say, $m<200$. This can be explained in the following way: for small $m$, updating approximate solutions does not affect the CPU time so much. Other operations such as matrix-vector multiplications are now the most time-consuming parts since the three methods required almost the same number of iterations, see Fig.~4.1. From Proposition 2.2 $(I)$ and Proposition 3.2 $(I)$ we know  that in this case the shifted QMR\_SYM method and the shifted QMR\_SYM($B$) method require real matrix-real vector multiplications. On the other hand, the shifted COCG method requires real matrix-complex vector multiplications.  Moreover, dot products and vector additions of the complex symmetric Lanczos process used in the shifted QMR\_SYM method and the shifted QMR\_SYM($B$) method can be done in real arithmetic. Hence, the two methods converged much faster than the shifted COCG method.\\ \\ 
3.2. Example 2\\ \\
The second problem comes from the electronic structure calculation of bulk fcc Cu with 1568 atoms in \cite{Shifted_COCG} and is given as follows: 
\begin{eqnarray*}
(\sigma_\ell I - H)\xb^{(\ell)} = \eb_1, \ \  \ell=1,2,\ldots,m,
\end{eqnarray*}
where $\sigma_\ell=-0.5+(\ell-1+i)/1000$, $H\in R^{14112\times 14112}$ is a symmetric matrix with 3924704 entries,  $\eb_1=(1,0,\ldots,0)^\TT$, and $m=1501$.

\begin{figure}[htbp] 
\begin{center}
\includegraphics[width=0.8\linewidth, height=0.5\linewidth]{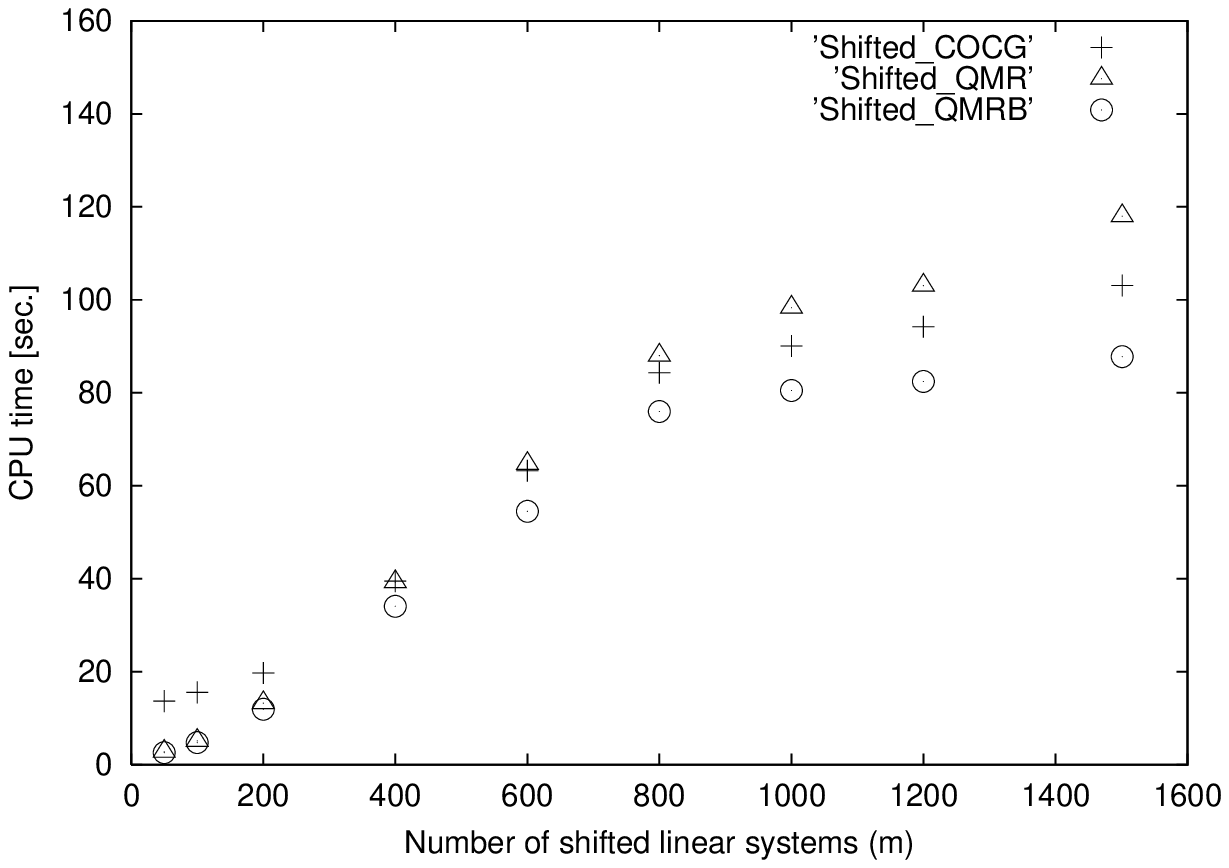}
\caption{Required CPU time given in seconds versus the number of shifted linear systems for each iterative method.} \vspace{2mm}
\includegraphics[width=0.8\linewidth, height=0.5\linewidth]{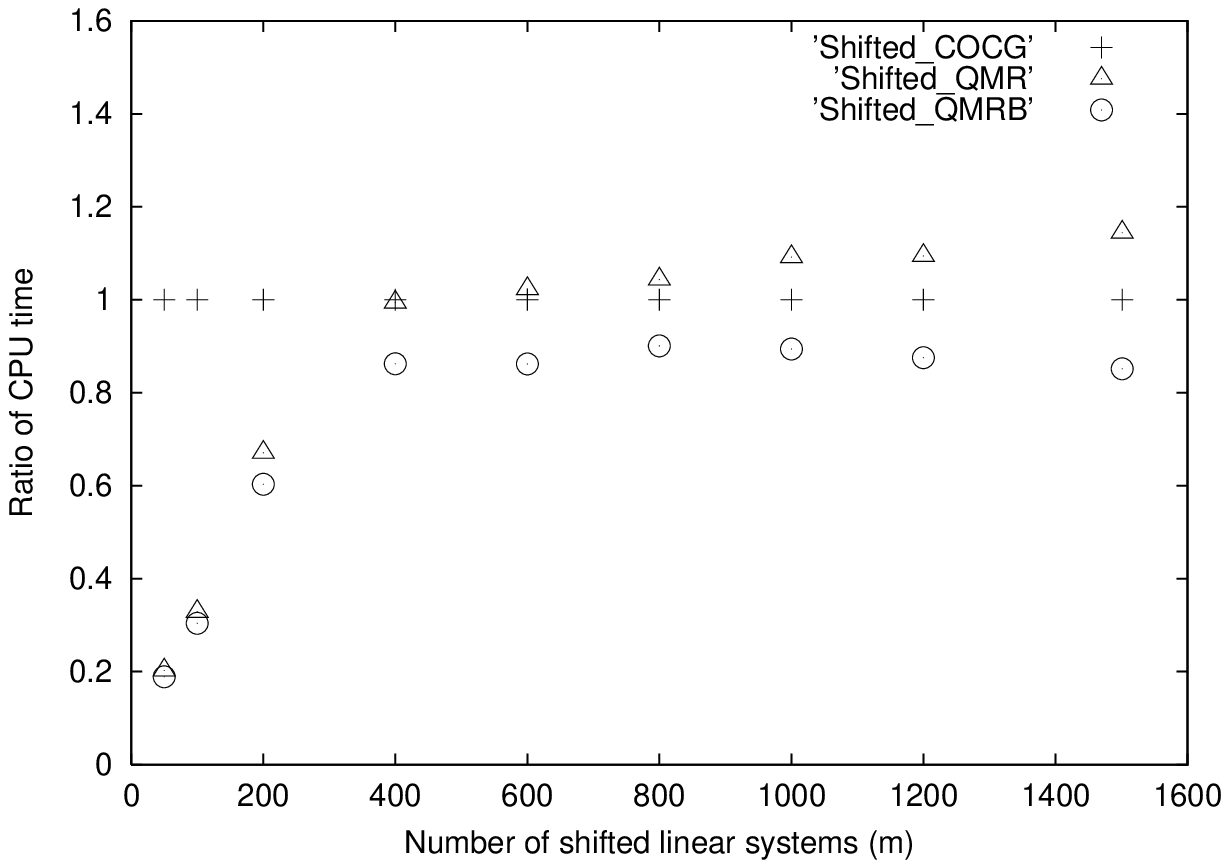}
\caption{The ratio of each computational time to the one of the shifted COCG method versus the number of shifted linear systems.} 
\end{center}
\end{figure}

Each computational time of the three methods for solving the $m$ shifted linear systems is shown in Fig.~4.5, where horizontal and vertical axes are the same as those used in Fig~4.3.  The ratio of each computational time to the computational time of the shifted COCG method is shown in Fig.~4.6. From Figs. 4.5, 4.6 we see that the size of this matrix is about 7 times larger than before, and the three methods show similar tendencies to those seen in the previous example.

\section{Concluding remarks}
In this paper, the shifted QMR\_SYM method was described as a specialization of the QMR method for general non-Hermitian shifted linear systems \cite{Freund}. The method has an advantage over the shifted COCG method in that there is no need to choose a suitable seed system.  On the other hand, we have found that for a large number of shifted linear systems, the most time-consuming part of the shifted QMR\_SYM method is updating approximate solutions, and this cost is higher than that of the shifted COCG method.  We therefore have proposed the weighted quasi-minimal residual approach with a specific weight for reducing the computational cost for updating approximate solutions. The resulting method, shifted QMR\_SYM($B$), also does not require to choose a suitable seed system, which is an advantage over the shifted COCG method. From numerical experiments we have learned that shifted QMR\_SYM and QMR\_SYM($B$) are competitive in comparison to the shifted COCG method. In particular, QMR\_SYM($B$) can be the method of choice for solving complex symmetric shifted linear systems with a large number of shifts that arise from large-scale electronic structure theory. In future work, numerical tests for general complex symmetric shifted linear systems will be done to evaluate the performance of the method.  \\

%

\section*{Acknowledgments}
We wish to express our gratitude to Roland W. Freund (UC Davis) for his fruitful comments during the conference at Harrachov in 2007. We are grateful to an anonymous referee for useful comments that substantially enhanced the quality of this paper. Finally, we would like to thank Martin H. Gutknecht (ETH Zurich) for careful reading of the manuscript.\\ \\

  \newpage

\end{document}